\documentclass[11pt]{article}   	
\pdfoutput=1 		

\usepackage{graphicx}				
\usepackage{amssymb}
\usepackage{amsmath}

\usepackage{subcaption}

\newtheorem{theorem}{Theorem}

\newcommand{\rr}{{\mathbf r}}
\newcommand{\s}{{\mathbf s}}

\newcommand{\nn}{{\mathbf n}}

\newcommand{\Reyuls}{{I{\relax\kern-.3em}R}}
\newcommand{\Phit}{\tilde{\Phi}}

\newcommand{\n}{\mathbf{n}}

\newcommand{\es}{\epsilon_s}
\newcommand{\ep}{\epsilon_p}
\newcommand{\ez}{\epsilon_0}

\title{A Hybrid Solver of Size Modified Poisson-Boltzmann Equation by Domain Decomposition, Finite Element, and Finite Difference}

\author{Jinyong Ying and Dexuan Xie\thanks{Corresponding author: dxie@uwm.edu (D. Xie)}}

\begin{document}
\maketitle

\begin{center}
 Department of Mathematical Sciences, University of Wisconsin-Milwaukee\\
Milwaukee, Wisconsin, USA, 53201-0413
\end{center}

\begin{abstract}
The size-modified Poisson-Boltzmann equation (SMPBE) is one important variant of the popular  dielectric model, the Poisson-Boltzmann equation (PBE), to reflect ionic size effects in the prediction of electrostatics for a biomolecule in an ionic solvent. In this paper, a new SMPBE hybrid solver is developed using solution decomposition, Schwartz's overlapped domain decomposition, finite element, and finite difference. It is then programmed as a software package in \texttt{C}, \texttt{Fortran}, and \texttt{Python} based on the state-of-the-art  finite element library  \texttt{DOLFIN} from the \texttt{FEniCS} project. This software package is well validated on a Born ball model with analytical solution and a dipole model with a known physical properties. Numerical results on six proteins with different net charges demonstrate its high performance. Finally, 
this new SMPBE hybrid solver is shown to be numerically stable and convergent in the calculation of electrostatic solvation free energy for 216 biomolecules and binding free energy for a DNA-drug complex.
\end{abstract}


Keywords: 
{\em Poisson-Boltzmann equation, finite element method, finite difference method, domain decomposition, electrostatic solvation free energy, binding free energy}


\section{Introduction}\label{intro}

The Poisson-Boltzmann equation (PBE) has been widely applied to the prediction of electrostatics for a biomolecule in an ionic solvent and the calculation of many biophysical quantities, such as electrostatic solvation and binding free energies \cite{baker2004poisson,chen1997monovalent,honig95}, due to the popularity of the PBE software APBS \cite{unni2011web},   DelPhi \cite{smith2012delphi}, PBSA\cite{luo2002accelerated,wang2013exploring}, UHBD \cite{uhbd91}, and PBEQ \cite{jo2008pbeq}. However, PBE has been known not to work properly in the prediction of ionic concentrations since it simply treats  each ion as a volumeless point. 
To reflect ionic size effects, one variant of PBE, called the size modified PBE (SMPBE), was proposed  based on the assumption that each ion and each water molecule occupy the same space of a cube with side length $\Lambda$ \cite{borukhov1997steric}. It was revisited in \cite{BoLi2009a} under variational principle to yield a slightly different definition, and shown to be optimal in the sense of minimizing a traditional electrostatic energy.  Another modification of SMPBE was given in \cite{chu2007evaluation}. We also noted that  there existed a nonuniform size modified PBE model \cite{BoLi2009a} and several other SMPBE models studied by the techniques of Monte Carlo and mean-field  \cite{andresen2004spatial,coalson1995statistical},  a generalized Poisson-Fermi distribution \cite{tresset2008generalized}, and statistical mechanics  \cite{boschitsch2012formulation}. 

Even so, the simple SMPBE models from \cite{borukhov1997steric,chu2007evaluation,BoLi2009a} remain to be the valuable ones for biomolecular applications due to their similarity to PBE in equation form. Currently, the SMPBE model given in \cite{chu2007evaluation} was solved by a finite difference scheme from APBS \cite{baker01}, which had problems of low accuracy and numerical instability\cite{Yi_Xie2014,zhou2008highly}. In \cite{chaudhry2011finite}, a finite element algorithm was proposed to solve the SMPBE defined in \cite{borukhov1997steric}, but its numerical tests were limited to a Born ion model and a small molecule with three atoms. As a generalization of our PBE finite element solver \cite{xiePBE2013}, we recently developed an effective SMPBE finite element solver \cite{Li_Xie2014b} according to the definition given in \cite{BoLi2009a}, and showed that it worked well for proteins with different net charges. The purpose of this paper is to modify it  as  a hybrid solver to further improve its performance. 

Our new SMPBE hybrid solver was motivated from the following observations. The finite element method with an unstructured interface-matched tetrahedral mesh can be effectively used to deal with the interface conditions of SMPBE defined on an interface $\Gamma$ with very irregular geometry. However, its implementation requires extra arrays to store mesh data and coefficient matrices of  finite element linear systems. Multigrid algorithms for solving each involved finite element linear system may be less efficient than the ones for solving a corresponding finite difference linear system \cite{trottenberg2000multigrid}. For example, we tested the preconditioned conjugate gradient method (PCG) with an algebraic multigrid preconditioner from the scientific computing library \texttt{PETSc}  \cite{petsc-user-ref}; it was found to take much more CPU run time than the PCG using incomplete LU preconditioning (PCG-ILU). This is the reason why PCG-ILU was selected as the default linear finite element solver in \cite{Li_Xie2014b}. On the other hand, a finite difference method using a uniform mesh can be solved very efficiently by a geometric multigrid scheme without storing any coefficient matrix or mesh data. However, developing a finite difference method for solving SMPBE can be very difficult despite some advances made in the treatment of the interface conditions in the case of PBE \cite{boschitsch2011fast,wang2013exploring}.  

We further noted that the calculation amount of the  SMPBE finite element  solver mainly came from solving a linear interface problem (see \eqref{Psi}) for $\Psi$ and a linear variational problem  (see \eqref{direction}) for a search direction, $p_{k}$, of a modified Newton variational minimization algorithm (see \eqref{Newton4Phit}) for computing $\Phit$. Here a sum of $\Psi$ and $\Phit$ with a known function $G$ (see \eqref{G-def}) gives a numerical solution $u$ of SMPBE. Hence, we only need to construct a hybrid algorithm for computing $\Psi$ and $p_{k}$ to modify  the SMPBE finite element solver into a hybrid solver. 

To do so, one key step is to reformulate  the linear interface problem  that defines the search direction $p_k$ (see \eqref{thm-p}) from a variational form into a differential form. In this paper, this reformulation is done in Theorem 4.1. We then construct two overlapped domain decomposition schemes for solving $\Psi$ (see \eqref{box-iterates-Psi}) and $p_k$ (see \eqref{box-iterates-P}), respectively,   based on a special overlapped seven box partition (See Subsection \ref{boxpartitionscheme}). In this box partition, the central box contains the protein region $D_{p}$, and is surrounded by six overlapped neighboring boxes. In order to simplify  the data exchange between any two neighboring boxes, we next construct a special mesh of the central box, which mixes an unstructured interface-matched tetrahedral mesh with a regular tetrahedral mesh. Furthermore, a finite element scheme is applied to the central box to solve a linear interface problem while a finite difference scheme is applied to each neighboring box to solve a Poisson (or Poisson-like) boundary value problem. In this way, we obtain the two finite element and finite difference hybrid algorithms, one for computing $\Psi$ and the other for $p_k$. Using them, we modify the SMPBE finite element solver into a new hybrid one. 

From the standard Schwartz's domain decomposition theory \cite{xu1998some} it can be known that our overlapped domain decomposition scheme has a fixed rate of convergence for a fixed over-relaxation parameter $\omega$. Thus, its performance mainly depends on the performance of a linear iterative scheme within each box. 
In our SMPBE hybrid solver, we retain the PCG-ILU as the finite element solver within the central box, and develop an ``optimal'' scheme, the PCG using multigrid V-cycle preconditioning  (PCG-MG), to solve each finite difference linear system within each neighboring box.
Because of our special interface-matched tetrahedral mesh of the central box, the data exchange between the finite element and finite difference methods can be carried out easily and efficiently. 
 
We programmed this new SMPBE hybrid scheme in \texttt{C}, \texttt{Fortran}, and \texttt{Python} as a modification of the finite element program package reported in \cite{Li_Xie2014b}. The new program parts include a \texttt{Fortran} program of PCG-MG and a special mesh generation program for the central box, which we developed based on our revised version of the molecular surface and volumetric mesh generation program package \texttt{GAMer} \cite{yu2008feature}. Although PCG-MG is a well known scheme, we did not find any software that is suitable for our case. Thus, we programed it ourselves. Our PCG-MG program was done based on the \texttt{BLAS} library ({\em http://www.netlib.org/blas/}). In this implementation, all the required memory arrays are preallocated. It does not require any memory array to store mesh data or  coefficient matrices of  finite difference linear systems. As a separate software, it can also be easily adopted to solve a general Poisson-like boundary value problems on a rectangular box. 

We validated this new hybrid program package using a new SMPBE test model for arbitrarily multiple charges artificially constructed based on a Poisson model from \cite{PBEHybridModel}, whose analytical solution is given. Numerical results from these tests also confirmed that both PCG-MG and our special overlapped domain decomposition scheme had convergence rates independent of the mesh size $h$, and were very efficient in terms of CPU time and memory usage. To demonstrate that SMPBE is a better model than PBE in the prediction of ionic concentrations, we constructed a more interesting dipole test model  than a commonly-used Born ion ball model. In this dipole model, the solute region $D_p$ consists of two overlapped balls with the same radius but two opposite central charges. Our numerical results on this dipole model  showed that SMPBE can much better capture physical features of ionic solvent than PBE (see Figures~\ref{vectorfield} and \ref{concentration_dipole}). Furthermore, we made numerical experiments on six proteins with different net charges in a range from $-35e_c$ to $+6e_c$ to compare the performance of the new hybrid solver with that of the  finite element solver. From the numerical results of Table~\ref{protein-case} it can be seen that the total CPU runtime of the SMPBE finite element solver was reduced sharply up to $76\%$ by our new SMPBE hybrid solver. For example, it took only about 18 seconds for our hybrid solver to find a numerical solution of SMPBE on a mesh with 537,953 mesh points on one 3.7 GHZ processor of our Mac Pro Workstation with 64 GB memory.

Prediction of electrostatic solvation free energy is one important application of SMPBE. It can also be used to validate a SMPBE solver since a reliable numerical solver shall produce a sequence of energy values toward a limit or around an experimental value from a chemical laboratory as the mesh size $h\to 0$. To check the numerical behavior of our new solver, we constructed six sets of meshes with the numbers of mesh points increasing from 28,166 to 2,086,780. We then calculated the electrostatic solvation free energies for 216 biomolecules we obtained from a protein
database maintained by Prof. Ray Luo on the homepage {\em http://rayl0.bio.uci.edu/rayl/}. These test results (see Figure~\ref{SolVConvergence}) confirm numerically that our SMPBE hybrid solver has good properties in numerical stability and convergence.

Finally, the calculation of binding free energy was done by our hybrid solver for a DNA-drug complex from \cite{breslauer1987enthalpy,fenley2010revisiting}. In these tests, we calculated binding free energies using six different sets of meshes, whose number of mesh points varied from about 85,000 to 4,250,000. The binding free energies and their scaled slopes with respect to the ionic strength $I_s$ were calculated for 11 different values of $I_{s}$ in a range from $e^{-3}$ to $e^{-1}$ (i.e., about 0.05 to 0.37). They were found to behave stably on different meshes and well matched the chemical experiment data from a chemical laboratory as given in \cite{breslauer1987enthalpy}.

The remaining parts of the paper are organized as follows. Section \ref{SMPBEReview} introduces SMPBE and its solution decomposition. Section \ref{FEP_Review} reviews the SMPBE finite element solver. Section \ref{Hybrid_Review} presents the new SMPBE hybrid solver. Finally,  the numerical results are reported in Section \ref{NumTests}. 

\section{SMPBE and its solution decomposition}
\label{SMPBEReview}

Let $D_{p}$ denote a solute region hosting a protein molecule (or other biomolecules such as DNA and RNA) with $n_{p}$ atoms,  $D_{s}$  a solvent region,  $\Gamma$ an interface between $D_{p}$ and $D_{s}$, and  $\Omega$ a sufficiently large domain  to satisfy that 
\[D_{p} \subset \Omega \quad \mbox{and}\quad  \Omega=D_p\cup D_s\cup \Gamma. \]
Both $D_{p}$ and $D_{s}$ are treated as continuum media with dielectric constants $\ep$ and $\es$, respectively.
To reflect ionic size effects, each ion and each water molecule are assumed to occupy the same volume of a cube with side length $\Lambda$. 

For a symmetric 1:1 ionic solvent (e,g., a salt solution with sodium (Na$^+$) and chloride (Cl$^-$) ions), the SMPBE model  has the following dimensionless form 
\begin{equation}
\label{SMPBE}
\left\{\begin{array}{ll}
-\ep\Delta u(\rr)=\alpha\displaystyle\sum_{j=1}^{n_p}z_j\delta_{\rr_j}, &\qquad \rr\in D_p, \\
-\es\Delta u(\rr) + \displaystyle\frac{\kappa^2\sinh(u)}{1+2M\Lambda^3\cosh(u)}=0 , &\qquad \rr\in D_s,\\
    u(\s^+)=u(\s^-), \quad
 \displaystyle \es\frac{\partial u(\s^+)}{\partial\nn(\s)}=\ep\frac{\partial u(\s^-)}{\partial\nn(\s)}, &\qquad  \s\in\Gamma,\\
 u(\s)=g(\s), &\qquad \s\in\partial\Omega,
\end{array}\right.
\end{equation}
where $\alpha$, $\kappa^2$, and $M$ are constants, $\rr_{j}$ and $z_{j}$ are the position and charge number of the $j$th atom of the protein, respectively, $g$ is a given boundary function, $\partial\Omega$ denotes the boundary of $\Omega$, $\delta_{\rr_j}$ is the Dirac delta distribution at point $\rr_{j}$, and $\nn(\s)$ is the unit outward normal vector of $D_{p}$.

Under the SI (Le Syst\`eme International d$^{'}$Unit\'es) units, for a domain $\Omega$ given in angstroms (\AA), the constants $\alpha$, $\kappa^2$, and $M$ are given by
\begin{equation}
\label{alpha-kappa}
 \alpha =  \frac{10^{10}e_{c}^{2}}{\ez k_{B}T}, \quad  \kappa^{2}=2I_{s} \frac{10^{-17}N_{A}e_{c}^{2}}{\ez k_{B}T},
 \quad\mbox{and}\quad  M=10^{-27} N_A I_s,
\end{equation}
where $k_B$, $T$, $N_A$, $I_{s}$, $e_c$ and $\ez$ denote the Boltzmann constant, the absolute temperature, the Avogadro number, the ionic strength in mole per liter, the electron charge, and the permittivity of vacuum, respectively.  For $T=298.15$ and $I_s=0.1$, we can get  
\begin{equation}\label{alpha-kappa2}
\alpha=7042.94, \quad \kappa^2=0.84827, \quad\text{and}\quad  M=6.0221\times 10^{-5},
\end{equation}
which will be used in our numerical tests in Section \ref{NumTests}. With these quantities, the corresponding solution $u$ is the electrostatic potential in units $k_BT/e_c$. 

To overcome the difficulties caused by the Dirac delta distributions $\delta_{\rr_{j}}$, a solution decomposition has been proposed in  \cite{Li_Xie2014b} to split the solution $u$ of \eqref{SMPBE} into the form  
\begin{equation}\label{solutionSplit}
u=G+\Psi+\Phit,
\end{equation}
where $G$ is given by the expression
\begin{equation}
\label{G-def}
    G(\rr)=\frac{\alpha}{4\pi\ep}\displaystyle \sum_{j=1}^{n_p}\frac{z_j}{|\rr-\rr_j|},
\end{equation}
$\Psi$ is a solution of the linear interface boundary value problem   
 \begin{equation}\label{Psi}
\left\{\begin{array}{cl}
 \Delta\Psi(\rr)=0, & \qquad\rr\in D_p\cup D_s, \\
\Psi(\s^+)=\Psi(\s^-),&\qquad\s\in\Gamma,\\
\displaystyle\es\frac{\partial\Psi(\s^+)}{\partial\nn(\s)}=\ep\frac{\partial\Psi(\s^-)}{\partial\nn(\s)}+(\ep-\es)\frac{\partial G(\s)}{\partial \nn(\s)}, & \qquad \s\in\Gamma, \\
\Psi(\s)=g(\s)-G(\s), &\qquad \s\in\partial\Omega,
\end{array}\right.
\end{equation}
and $\Phit$ is a solution of the nonlinear interface boundary value problem 
\begin{equation}
\label{Phit}
\left\{\begin{array}{cl}
\Delta\Phit(\rr)=0, &\qquad\rr\in D_p,\\
-\es\Delta \Phit(\rr)+\displaystyle\frac{\kappa^2\sinh(G+\Psi+\Phit)}{1+2M\Lambda^3\cosh(G+\Psi+\Phit)}=0 , &\qquad \rr\in D_s,\\
 \displaystyle\Phit(\s^+)=\Phit(\s^-), \quad  \es\frac{\partial\Phit(\s^+)}{\partial\nn(\s)}=\ep\frac{\partial\Phit(\s^-)}{\partial\nn(\s)}, &  \qquad\s\in\Gamma,\\
\Phit(\s)= 0,  &\qquad \s\in \partial\Omega.
\end{array}\right.
\end{equation}
Here $\displaystyle\frac{\partial G(\s)}{\partial \nn(\s)}=\nabla G\cdot \n$ with $\nabla G$ being given by 
\begin{equation}
\label{PG-def}
     \nabla G(\rr)= -\frac{\alpha}{4\pi\ep}\displaystyle \sum_{j=1}^{n_p}z_j\frac{\rr-\rr_j}{|\rr-\rr_j|^3}.
\end{equation}

\section{The SMPBE finite element solver}
\label{FEP_Review}

In this section, we briefly review the SMPBE finite element solver proposed in \cite{Li_Xie2014b}. 

Let $\mathcal{M}$ denote a finite element function space. It is assumed to be a subspace of the usual Sobolev function space $H^1(\Omega)$. A finite element solution $\Psi$ of \eqref{Psi} is defined by the following linear variational problem:\vspace{0.8em}

\centerline{Find $\Psi\in  {\mathcal M}$ with $\Psi|_{\partial\Omega}=g-G$ such that }
\vspace{-1.0em}
\begin{equation}
\label{PsiW}
a(\Psi, v)=(\ep-\es)\int_{D_s}\nabla G(\rr)  \cdot \nabla v(\rr) d\rr  \quad \forall v \in {\mathcal M}_{0},
\end{equation}
where $\mathcal{M}_{0}$ is a subspace of $\mathcal{M}$ defined by
 $$\mathcal{M}_{0}=\{v\in \mathcal{M} \; | \; v=0 \mbox{ on } \partial \Omega\},$$
and $a(u,v)$ is a bilinear form defined by
\begin{equation}\nonumber
a(u, v) = \ep \int_{D_p}\nabla u(\rr)  \cdot \nabla v(\rr) d\rr +\es \int_{D_s}\nabla u(\rr) \cdot \nabla v(\rr) d\rr.
\end{equation}
Clearly, $\mathcal{M}_{0}$ is a subspace of the Sobolev function space $H^1_0(\Omega)$. 

After $\Psi$ is found, we treat $U=\Psi + G$ as a given function. The nonlinear interface problem \eqref{Phit} can then be formulated to the nonlinear variational minimization problem: 
$ \displaystyle\min_{v\in\mathcal{M}_0}J(v),$
where $J$ is given by
$$J(v) = \frac{1}{2}a(v,v)+\frac{\kappa^{2}}{2M\Lambda^3}\int_{D_s}\ln(1+2M\Lambda^3\cosh(U+v))d\rr.$$
The above minimization problem is solved by the modified Newton minimization method
\begin{equation}
\label{Newton4Phit}
   \Phit^{(k+1)}=\Phit^{(k)}+\lambda_kp_k, \quad k=0,1,2,\cdots,
\end{equation}
where $\Phit^{(0)}$ is an initial guess, $\lambda_k$ is a steplength selected by using a line search algorithm to satisfy the following condition
\begin{equation}\label{DecentCondition}
J(\Phit^{(k)}+\lambda_kp_k)\leq J(\Phit^{(k)}) \quad\mbox{or}\quad  \|J^{\prime}(\Phit^{(k)}+\lambda_kp_k)\|\leq \|J^{\prime}(\Phit^{(k)})\|,
\end{equation}
and $p_{k}$ is a search direction satisfying the Newton equation in the variational form
\begin{equation}
\label{direction}
J''(\Phit^{(k)})(p_k, v) = -J'(\Phit^{(k)})v  \quad \forall v\in \mathcal{M}_{0},
\end{equation}
where $J'(\Phit)$ and $J''(\Phit)$ are the first and second Fr\'echet-derivative of $J$ at $\Phit$, respectively, which are defined by
$$J'(\Phit)v = a(\Phit,v)+\kappa^2\int_{D_s}\frac{\sinh(U+\Phit)}{1+2M\Lambda^3\cosh(U+\Phit)}vd\rr \quad \forall v\in H^1_0(\Omega),$$
$$J''(\Phit)(p,v) = a(p,v)+\kappa^2\int_{D_s}\frac{2M\Lambda^3+\cosh(U+\Phit)}{(1+2M\Lambda^3\cosh(U+\Phit))^2}p vd\rr  \hspace{1em}\forall p,v\in H^1_0(\Omega).$$
In computation, the Newton equation \eqref{direction} is solved by PCG-ILU such that the relative or the absolute residual norm is less than $10^{-8}$ by default. 

When the modified Newton iterate $\Phit^{(k)}$ reaches the iteration stop rule  
 \[  \|J^{\prime}(\Phit^{(k)})\|<\epsilon  \quad  \mbox{ (by default, $\epsilon=10^{-7}$)},\]
$\Phit^{(k)}$ is output as an finite element solution $\Phit$ of \eqref{Phit} on $\mathcal{M}_{0}$. A finite element solution $u$ of SMPBE is then constructed by the solution decomposition formula $u=G+\Psi+\Phit$ on ${\cal M}$.

\section{The new SMPBE hybrid solver}
\label{Hybrid_Review}

In this section, we present a new finite element and finite difference hybrid algorithm for solving SMPBE. 
One key step to construct a SMPBE hybrid scheme is to reformulate the Newton equation \eqref{direction} from the variational form into a linear interface problem in strong differential sense. We first complete this reformulation in Theorem~\ref{thm-p}. We then describe the new hybrid solver in four subsections for clarity.

\vspace{2mm}
\begin{theorem}
\label{thm-p}
Let  $w=\Phit^{(k)}+\Psi + G$ be given. If  $p$, $\Phit^{(k)} \in H^1_0(\Omega)\cap C^2(D_p)\cap C^2(D_s)$, then the Newton equation \eqref{direction} is equivalent to the following interface problem
\begin{subequations}
\label{SMPBE_search}
\begin{align}
\label{SMPBE_search1}
&-\Delta p(\rr)= \Delta \Phit^{(k)}(\rr), \hspace{4.3cm}\qquad \rr\in D_p,  \\
\label{SMPBE_search2}
&-\es\Delta p(\rr) + \kappa^2\frac{2M\Lambda^3+\cosh(w)}{(1+2M\Lambda^3\cosh(w))^2}p(\rr)= f_{s}(\rr),\quad  \rr\in D_s,\\
 \label{SMPBE_search3}
&\displaystyle    p(\s^+)=p(\s^-),\quad 
\displaystyle \es\frac{\partial p(\s^+)}{\partial\nn(\s)}-\ep\frac{\partial p(\s^-)}{\partial\nn(\s)}=g_{\Gamma}(\s), \hspace{2.5mm}\quad \s\in\Gamma,\\
\label{SMPBE_search4}
&p(\s)=0, \hspace{6.35cm}\qquad \s\in\partial\Omega,
\end{align}
\end{subequations}
where  $f_{s}$ and $g_{\Gamma}$ are given by  
\begin{equation}
\label{fg-def}
  f_{s}(\rr) =
\es \Delta \Phit^{(k)} -\frac{\kappa^2 \sinh(w)}{1+2M\Lambda^3\cosh(w)}, \quad
   g_{\Gamma}(\s)=\ep\frac{\partial \Phit^{(k)}(\s^-)}{\partial\nn(\s)} -\es\frac{\partial \Phit^{(k)}(\s^+)}{\partial\nn(\s)}.
\end{equation}
\end{theorem}

\vspace{2mm}

{\em Proof.}
 We only show the derivation of \eqref{SMPBE_search} from the variational form \eqref{direction} since the proof of the converse is easy. For any $v \in H^1_0(\Omega)$ satisfying $v=0$ on $D_{s}$ and $v \in \mathbb{C}^{\infty}_0(D_p)$, from \eqref{direction} we can get  
\[   \int_{D_p}(\Delta p+\Delta \Phit^{(k)})v d\rr= 0 \quad \forall v\in \mathbb{C}^{\infty}_0(D_p),\]
from which it implies equation \eqref{SMPBE_search1}.

Next, for any $v \in H^1_0(\Omega)$ satisfying $v=0$ on $D_{p}$ and $v \in \mathbb{C}^{\infty}_0(D_s)$, \eqref{direction} can be reduced to the form
\begin{align*}
&\es\int_{D_s}\nabla p\cdot \nabla v d\rr+\kappa^2\int_{D_s}\frac{2M\Lambda^3+\cosh(w)}{(1+2M\Lambda^3\cosh(w))^2}p vd\rr\\
=&-\es\int_{D_s}\nabla \Phit^{(k)}\cdot \nabla v d\rr -\kappa^2\int_{D_s}\frac{\sinh(w)}{1+2M\Lambda^3\cosh(w)}vd\rr.
\end{align*}
By Green's identity, the above equality can be reformulated as
\begin{align*}
& \int_{D_s} \left[-\es\Delta p+\kappa^2 \frac{2M\Lambda^3+\cosh(w)}{(1+2M\Lambda^3\cosh(w))^2}p \right]vd\rr \\
= & \es\int_{D_s}\Delta \Phit^{(k)} vd\rr  -\kappa^2\int_{D_s}\frac{\sinh(w)}{1+2M\Lambda^3\cosh(w)}vd\rr
 =\int_{D_s}f_{s}(\rr)v(\rr)d\rr,
\end{align*}
from which we can obtain equation \eqref{SMPBE_search2}.

Furthermore, applying Green's identity to the two terms of $a(\cdot,\cdot)$ for $v\in H^1_0(\Omega)$ in $D_p$ and $D_s$, respectively, we can reformulate \eqref{direction} as 
\begin{align*}
&\int_{\Gamma}\left[ \ep \frac{\partial p(\s^-)}{\partial \nn}-\es\frac{\partial p(\s^+)}{\partial \nn}\right] vd\s -\ep\int_{D_p}\Delta p vd\rr  - \int_{D_s} \left[ \es\Delta p - \kappa^2 \frac{2M\Lambda^3+\cosh(w)}{(1+2M\Lambda^3\cosh(w))^2}p\right] vd\rr \\
& = \int_{\Gamma}\left[\es\frac{\partial \Phit^{(k)}(\s^+)}{\partial \nn}-\ep\frac{\partial \Phit^{(k)}(\s^-)}{\partial \nn}\right]vd\s +\ep\int_{D_p}\Delta \Phit^{(k)} vd\rr \\
&\quad + \int_{D_s}\left[\es \Delta \Phit^{(k)}-\kappa^2\frac{\sinh(w)}{1+2M\Lambda^3\cosh(w)}\right]vd\rr.
\end{align*} 
Applying \eqref{SMPBE_search1} and \eqref{SMPBE_search2} to the above identity leads to the interface condition \eqref{SMPBE_search3}. The boundary condition \eqref{SMPBE_search4} is natural because $p\in H^1_0(\Omega)$.
This completes the proof.

\vspace{2mm}

Because of Theorem~\ref{thm-p}, both $\Psi$ and the search direction $p_k$ can now be calculated by a finite difference method. This makes it possible for us to modify the SMPBE finite element solver as a finite difference solver. But, to combine the advantages of finite element and finite difference solvers together, in the remaining part of this section, we present a novel finite element and finite difference hybrid algorithm for solving SMPBE.

\subsection{A special overlapped box partition}
\label{boxpartitionscheme}
In this subsection, we  describe a scheme for constructing the domain $\Omega$ and a special overlapped box partition of $\Omega$. Here, a protein region $D_{p}$ is given, and a box $D$ is selected to satisfy $D_{p}\subset D$.

For simplicity, we set $D=\prod_{i=1}^{3}(a_{i},b_{i})$ as a cube with side length $b_i-a_i=L$ for $i=1,2,3$.
We then define the mesh size $h$ and the  two parameters $\tau$ and $\eta$ by
\begin{equation}
\label{nmnu}
 h=L/2^{n}, \quad \tau = 2^{m} h, \quad \eta = \mu L/2,
\end{equation}
where $n$, $m$, and $\mu$ are three positive integers to be selected according to the need of calculation, and $\tau<\eta$. 
We next construct $\Omega$ as a cubic domain and its seven overlapped subdomains, $\Omega_{i}$, $i=1,2,\cdots,7$, as follows:   
\begin{equation}
\label{domain_construct}
\begin{split}
\Omega &= \prod_{i=1}^{3}(a_{i}-\eta,b_{i}+\eta), \hspace{1em} \Omega_{7} = \prod_{i=1}^{3}(a_{i}-\tau,b_{i}+\tau),\\
\Omega_1&=(a_1-\eta,b_1+\eta)\times(a_2-\eta,b_2+\eta)\times(a_3-\eta,a_3),\\
\Omega_2&=(a_1-\eta,b_1+\eta)\times(a_2-\eta,a_2)\times(a_3-\tau,b_3+\tau),\\
\Omega_3&=(a_1-\eta,a_1)\times(a_2-\tau,b_2+\tau)\times(a_3-\tau,b_3+\tau),\\
\Omega_4&=(b_1,b_1+\eta)\times(a_2-\eta,b_2+\tau)\times(a_3-\tau,b_3+\tau),\\
\Omega_5&=(a_1-\eta,b_1+\eta)\times(b_2,b_2+\eta)\times(a_3-\tau,b_3+\tau),\\
\Omega_6&=(a_1-\eta,b_1+\eta)\times(a_2-\eta,b_2+\eta)\times(b_3,b_3+\eta).
\end{split}
\end{equation}
Here, $\Omega_7$ is the central box, which is surrounded by the other 6 boxes. Clearly, we have
\begin{equation}
\label{box-conditions}
 (a) \quad D_p\subset D\subset\Omega_7;  \qquad  (b)\quad  \Omega\setminus D = \bigcup_{j=1}^{6}\Omega_j; \qquad (c)\quad \Omega=\bigcup_{j=1}^{7}\Omega_j.
\end{equation}
The positions and ordering numbers of these seven boxes are referred to \cite{ying2015new}.
        
\subsection{Two overlapped box  iterative methods}
Based on the above special overlapped box partition, we construct an overlapped box iterative method for computing $\Psi$ as follows: \\
\indent{For $m=1$, $2$, $3$, $\cdots$,  }
\begin{equation}
\label{box-iterates-Psi}
\Psi^{(m)}_{i} = (1-\omega) \Psi_{i}^{(m-1)}+\omega \widehat{\Psi}_{i} \quad \mbox{ on $\; \Omega_{i}\;$ for  } i=1,2,\cdots,7,
\end{equation}
where  
 $\Psi^{(0)}_{i}$ is an initial iterate, $\omega\in(1,2)$ is the over-relaxation parameter,  $\widehat{\Psi}_{i}$ with $i=1$ to 6 denotes a solution of the Poisson boundary value problem: 
\begin{equation}
\label{PDE_fd1}
\left\{\begin{array}{ll}
 \Delta \Psi(\rr) = 0  &\quad\mbox{in }\Omega_{i}, \\
\Psi(\s)=\Psi_{j}^{(m-1)}(\s) &\quad \mbox{on $\partial \Omega_{i}\cap\Omega_{j}$ if }\partial \Omega_{i}\cap\Omega_{j} \neq \varnothing \mbox{ for  $j=i+1$ to 7},\\
\Psi(\s)=\Psi_{j}^{(m)}(\s) &\quad \mbox{on $\partial \Omega_{i}\cap\Omega_{j}$ if }\partial \Omega_{i}\cap\Omega_{j} \neq \varnothing \mbox{ for   $j=1$ to $i-1$},\\
\Psi(\s)=g(\s)-G(\s) &\quad \mbox{on } \partial \Omega_{i}\cap\partial\Omega,
\end{array}\right.
\end{equation}
and $\widehat{\Psi}_{7}$ is a solution of the linear interface problem:
\begin{equation}
\label{PDE_fe1}
\left\{\begin{array}{ll}
 \Delta \Psi(\rr) =0 \quad\quad\quad  \mbox{ in }  \Omega_7 \setminus \Gamma,& \\
     \Psi(\s^+)=\Psi(\s^-), \quad
   \es\frac{\partial \Psi(\s^+)}{\partial\nn(\s)}=\ep\frac{\partial \Psi(\s^-)}{\partial\nn(\s)}+(\ep-\es)\frac{\partial G(\s)}{\partial \nn(\s)} &\mbox{ on } \Gamma,\\
\Psi(\s)=\Psi_{j}^{(m)}(\s) \quad \mbox{ on }\partial\Omega_7\cap\Omega_{j} \mbox{ for $j=1$ to 6}.&\\
\end{array}\right.
\end{equation}

After $\Psi$ and $\Phit^{(k)}$ are computed, we similarly construct another overlapped box iterative method for computing the search direction $p_{k}$ of \eqref{Newton4Phit} as follows: 
 For $m=1,2,3, \cdots$,  
\begin{equation}
\label{box-iterates-P}
p^{(m)}_{i} = (1-\omega) p_{i}^{(m-1)}+\omega \widehat{p}_{i} \quad \mbox{ on $\; \Omega_{i}\;$ for  } i=1,2,\cdots,7,
\end{equation}
where 
 $p^{(0)}_{i}$ is an initial iterate, $\omega\in(1,2)$ is the over-relaxation parameter,  $\widehat{p}_{i}$ with $i=1$ to 6 is a solution of the Poisson-like boundary value problem: 
\begin{equation}
\label{PDE_fd2}
\left\{\begin{array}{ll}
-\es\Delta p(\rr) + \displaystyle\kappa^2\frac{2M\Lambda^3+\cosh(\Phit^{(k)}+\Psi + G)}{(1+2M\Lambda^3\cosh(\Phit^{(k)}+\Psi + G))^2}p(\rr)= f_{s}(\rr) \quad\mbox{in }\Omega_{i}, \\
p(\s)=p_{j}^{(m-1)}(\s) \quad \mbox{on $\partial \Omega_{i}\cap\Omega_{j}$ if }\partial \Omega_{i}\cap\Omega_{j} \neq \varnothing  \mbox{ for  $j=i+1$ to 7},\\
p(\s)=p_{j}^{(m)}(\s) \quad \quad \mbox{on $\partial \Omega_{i}\cap\Omega_{j}$ if }\partial \Omega_{i}\cap\Omega_{j} \neq \varnothing \mbox{ for   $j=1$ to $i-1$},\\
p(\s)=0 \qquad\quad \hspace{1.15em}  \mbox{ on } \partial \Omega_{i}\cap\partial\Omega,
\end{array}\right.
\end{equation}
and $\widehat{p}_{7}$ is a solution of the linear interface problem:
\begin{equation}
\label{PDE_fe2}
\left\{\begin{array}{ll}
 -\Delta p(\rr) =\Delta \Phit^{(k)}(\rr) \quad \mbox{ in }  D_{p},& \\
 -\es\Delta p(\rr) +\displaystyle \kappa^2\frac{2M\Lambda^3+\cosh(\Phit^{(k)}+\Psi + G)}{(1+2M\Lambda^3\cosh(\Phit^{(k)}+\Psi + G))^2}p(\rr)= f_{s}(\rr)
 & \quad  \mbox{ in }   D_s\cap \Omega_7,\\
     p(\s^+)=p(\s^-), \quad
   \es\frac{\partial p(\s^+)}{\partial\nn(\s)}=\ep\frac{\partial p(\s^-)}{\partial\nn(\s)}+g_{\Gamma}(\s) &\quad \mbox{ on } \Gamma,\\
p(\s)=p_{j}^{(m)}(\s) \quad \mbox{ on }\partial\Omega_7\cap\Omega_{j} \mbox{ for $j=1$ to 6},&\\
\end{array}\right.
\end{equation}
where $f_{s}$ and $g_{\Gamma}$ have been given in \eqref{fg-def}.

Because the central box has been ordered as the last box $\Omega_{7}$, from \eqref{PDE_fe1} and \eqref{PDE_fe2} it can be seen that the updated values from the six neighboring boxes can be used to update the boundary value function. This may provide the interface problem on the central box with a better boundary value problem.

Due to the above two box iterative methods, the work amount of solving SMPBE mainly comes from the central box $\Omega_{7}$ and each neighboring boxe $\Omega_{i}$ for $i=1$ to 6 to solve a linear interface problem and a Poisson (or Poisson-like) boundary value problem, respectively. Clearly, different numerical methods can be applied to different boxes. This makes it possible for us to construct a more efficient SMPBE numerical solver than the SMPBE finite element solver.

\subsection{The new SMPBE hybrid solver}
Specifically, we use a finite element method to discretize each interface problem on the central box $\Omega_{7}$,  and a finite difference method to discretize each boundary value problem on each neighboring box $\Omega_{i}$ for $i=1$ to 6.
We then solve each finite element linear system by PCG-ILU and each finite difference linear system by PCG-MG. Consequently, we modify the SMPBE finite element solver as a new SMPBE hybrid solver. 
 For clarity,  this new hybrid solver is presented in Algorithm~1.
\vspace{2mm}

{\bf Algorithm 1} (The new  SMPBE hybrid solver). {\em Let an overlapped box partition of $\Omega$ be given as described in Subsection \ref{boxpartitionscheme}. A solution $u$ of the SMPBE model \eqref{SMPBE} is calculated approximately in the following five steps:
\begin{enumerate}
\item[Step 1.] Construct an interface-matched tetrahedral mesh for the central box $\Omega_{7}$ and a uniform mesh for each neighboring box $\Omega_{i}$ for $i\neq 7$ with a mesh size $h>0$.
\item[Step 2.]  Calculate $G$ on each box and $\nabla G$ on $\Omega_7$ according to \eqref{G-def} and  \eqref{PG-def}, respectively.
\item[Step 3.]  Calculate $\Psi$ by the overlapped box iterative method \eqref{box-iterates-Psi}. Here
\eqref{PDE_fd1} and \eqref{PDE_fe1} are approximated as finite difference and finite element linear systems and  solved by PCG-MG and  PCG-ILU, respectively, until the relative residual norm less than $10^{-8}$.
\item[Step 4.]   Calculate $\Phit$ by the modified Newton method in the following steps:
   \begin{enumerate} 
    \item Set $k=0$ and $\Phit^{(0)}=0$ (by default).
    \item Calculate the search direction $p_{k}$ by the overlapped box iterative method \eqref{box-iterates-P}. Here \eqref{PDE_fd2} and \eqref{PDE_fe2} are approximated as finite difference and finite element linear systems and  solved by PCG-MG and PCG-ILU, respectively, until the relative residual norm less than $10^{-8}$.
    \item Find the steplength $\lambda_{k}$ by a line search algorithm (starting with $\lambda_{k}=1$). 
     \item Define the modified Newton iterate $\Phit^{(k+1)}$ by $\Phit^{(k+1)}=\Phit^{(k)}+\lambda_kp_k.$
      \item Check the convergence: If  $ \|\Phit^{(k+1)}-\Phit^{(k)}\|\leq 10^{-7}$ (by default),
then $\Phit^{(k+1)}$ is set as  a solution $\Phit$ of the nonlinear interface problem \eqref{Phit}; otherwise, increase $k$ by 1 and go back to Step (b).
    \end{enumerate}
\item[Step 5.] Construct $u$ by the solution decomposition $u = G + \Psi + \Phit$.
\end{enumerate}
}

\subsection{The SMPBE hybrid solver program package.}
We programmed Algorithm 1 in \texttt{C}, \texttt{Fortran}, and \texttt{Python} as a software package. Similar to the SMPBE finite element solver package, the main program of our software was written in \texttt{Python} based on
the state-of-the-art  finite element library  \texttt{DOLFIN} from the \texttt{FEniCS} project \cite{logg2012automated,logg2010dolfin}. Each finite element equation is produced by using \texttt{DOLFIN}, and solved by PCG-ILU from the \texttt{PETSc} library \cite{petsc-user-ref}. The input file of the program is a PQR file of a protein molecule, which contains the positions $\rr_{j}$, charge numbers $z_{j}$, and radii of atoms as well as the related hydrogen atoms. The PQR file can be produced by the program tool \texttt{PDB2PQR} \cite{dolinsky2004pdb2pqr} from a PDB file of the protein, which can be downloaded from the Protein Data Bank (PDB) ({\em http://www.rcsb.org/}). 
This new software contains a mesh generation program we wrote in \texttt{C} to produce a special  mesh of $\Omega_7$, which mixes an unstructured interface-matched tetrahedral mesh of $D$ with a uniform tetrahedral mesh on the remaining part $\Omega_{7}\setminus D$. In this mesh program, the unstructured mesh is generated by our revised version of the molecular surface and volumetric mesh generation program package \texttt{GAMer} \cite{yu2008feature}, and the uniform mesh has the same mesh size $h$ as the one used in the construction of a uniform mesh of each neighboring box. This special mesh makes the data exchanges between two neighboring boxes easy and fast. 

A detailed description of PCG-MG can be found in \cite{trottenberg2000multigrid}. In this paper,
we programmed PCG-MG in \texttt{Fortran} without storing any mesh data or coefficient matrix of a finite difference linear system. In this program, all required memory arrays by PCG-MG are pre-allocated to improve performance. The multigrid precontitioner is defined by one V-cycle iteration with an initial guess of zero. Its main components are set as follows:
\begin{itemize}
\item The pre and post smoothing steps are defined by one forward and one backward Gauss-Seidel iteration, respectively.
\item The coarse grid meshes are  generated by using mesh sizes $h_{k}=2^{k-1}h$ for $k=1,2,\ldots,l$ satisfying  
$h=h_{1}<h_{2}<\ldots < h_{l}.$
Here $l$ is the number of mesh levels. A coarse grid mesh is set as the $l$-th grid mesh (or  the coarsest grid mesh) whenever it has only one or two interior mesh points in one or two coordinate directions.
\item The restriction matrix $I_{k-1}^{k}$ from a fine grid with mesh size $h_{k-1}$ to a coarse grid with mesh size $h_{k}$ and the prolongation matrix $I_{k}^{k-1}$ from the coarse mesh to the fine mesh are constructed by the standard full weight and trilinear interpolation techniques such that $I^{k-1}_{k} = 8 (I_{k-1}^{k})^T$ \cite[page 72]{trottenberg2000multigrid}.
\item The linear system on the coarsest grid mesh is solved by the LU factorization method since the size of each linear system on the coarsest grid mesh is small. 
\end{itemize}

For example, for $n=4$, $m=2$ and $\mu=4$,  we have $l=5$ for $\Omega_{1}$ and $\Omega_{6}$, and $l=4$ for the other four neighboring boxes. The linear systems on the coarsest grid mesh have only 16 unknowns for $\Omega_1$ and $\Omega_6$,
54 unknowns for $\Omega_2$ and $\Omega_5$, and 12 unknowns for $\Omega_3$ and $\Omega_4$.

We further wrote a \texttt{Fortran} subroutine for calculating the values of $G$ and $\nabla G$ at each mesh point in order to speed up the calculation of function $G$ and its gradient vector $\nabla G$. 
Our \texttt{Fortran} subroutines and   \texttt{C} programs were converted to Python external modules by the Fortran-to-Python interface generator \texttt{f2py} ({\em  http://cens.ioc.ee/projects/f2py2e/}) and \texttt{SWIG} ({\em http://www.swig.org}), respectively. Hence, they can be applied to the Python main program of our SMPBE hybrid software.

\section{Numerical results}
\label{NumTests}
In this section, we made numerical experiments to validate our new SMPBE program package and demonstrate its performance and applications. 
For simplicity, all the numerical tests were done by using $\ep=2.0$, $\es=80.0$, $\Lambda=3.11$, the values of three constants $\alpha$, $\kappa^2$, and $M$  given in \eqref{alpha-kappa2}, and other default parameter values for PCG-ILU, PCG-MG, and the modified Newton minimization algorithm. They were implemented on one processor of our Mac Pro Workstation with the 3.7 GHZ Quad-Core Intel Xeon E5 and 64 GB main memory.

\subsection{Validation tests}
In \cite{PBEHybridModel}, we have obtained the analytical solution 
$U$  of the Poisson equation, in a simple series expression,  for a spherical solute region containing multiple point charges in terms of Legendre polynomials. Using this analytical solution $U$, we construct a SMPBE test model as follows: 
\begin{equation}\label{SMBorn}
\left\{\begin{array}{ll}
-\epsilon_p\Delta u(\rr)=\alpha \displaystyle\sum_{n=1}^{n_p}  \delta_{\rr_j} &  \mbox{ in } D_p,\\
-\epsilon_s\Delta u(\rr)+\displaystyle\frac{\kappa^2\sinh(u)}{1+2M\Lambda^3\cosh(u)}=F_s(\rr) &\mbox{ in }  D_s,\\
\displaystyle u(\s^+)=u(\s^-),\quad \epsilon_s \frac{\partial u(\s^+)}{\partial \bold{n}}=\epsilon_p \frac{\partial u(\s^-)}{\partial \bold{n}} &\mbox{ on } \Gamma,\\
\displaystyle u(\rr)=U(\rr)&\mbox{ on  } \partial \Omega,
\end{array}\right.
\end{equation}
where  $D_p=\{\rr \; | \; |\rr|<a\}$ with $a>0$, $\Gamma=\{\rr \; | \;  |\rr|=a\}$, $\Omega$ is a box such that $D_s=\Omega-D_p-\Gamma$ is nonempty, and $F_s(\rr) =\kappa^2\sinh(U(\rr))[1+2M\Lambda^3\cosh(U(\rr))]^{-1}$, which can be understood as an excess charge density function. 
Clearly, $U(\rr)$ is also the analytical solution of the SMPBE test model \eqref{SMBorn}. 

In the numerical tests, we set radius $a=1$, $\Omega=(-6,6)^3$, and constructed an overlapped box partition of $\Omega$ using  $D=(-2,2)^3$, $\tau=1$, and $\eta=4$, which gave $\Omega_7=(-3,3)^3$. The over-relaxation parameter $\omega$ of our overlapped box iterative method was set as 1.275 and 1.225 in solving \eqref{Psi} for $\Psi$ and \eqref{SMPBE_search} for $p_{k}$, respectively. Three nested meshes with the mesh sizes $h=0.25$, $h/2$ (0.125), and $h/4$ (0.0625) were constructed for testing the convergence behavior of our hybrid solver. Their numbers of mesh points were 120887, 940247, and 7412989, respectively, including the numbers of mesh points from the finite element meshes of the central box $\Omega_7$, which were 18863, 145223, and 1136605. 

Table~\ref{nonlinearV} reports the numerical results for two validation test cases. In the first case, the unit ball region $D_{p}$ has only one central charge ($n_{p}=1$). Such a test model is often referred to as a Born ball test model. Its analytical solution $u$ is given by
\begin{equation}
\label{ASolution}
\displaystyle u(\rr)=\left\{\begin{array}{ll}
\displaystyle\frac{\alpha  }{4\pi}(\frac{1}{\epsilon_s}-\frac{1}{\epsilon_p})+\frac{\alpha  }{4\pi \epsilon_p|\rr|} & \mbox{ in } D_p,\\
\displaystyle\frac{\alpha  }{4\pi \epsilon_s|\rr|} & \mbox{ in } D_s.
\end{array}\right.
\end{equation}
 In the second case ($n_{p}=488$), we assigned the 488 atomic charges of a protein molecule (PDB ID: 2LZX) to the unit ball region $D_{p}$ through dividing each  atomic position $\rr_{j}$ by 19. Here the numbers of iterations for PCG-MG, PCG-ILU, and the overlapped box iterative methods were their averages over the total number of the linear systems solved.

\begin{table}[h]
\caption{Performance of our SMPBE hybrid  solver for the SMPBE test model \eqref{SMBorn} in relative solution errors and  average  iteration numbers (Iter.). 
}
\label{nonlinearV}
\centering
{
\begin{tabular}{|c||c|c|c|c|c|c|}
\hline
Mesh &   Error	 & PCG-MG Iter. & PCG-ILU &  Hybrid Box    & Newton  &  \\
size $h$ &  $\frac{\|u-u_h\|_{l^2(\Omega)}}{\|u\|_{l^2{\Omega}}}$	  &  on $\Omega_{i}$ ($i=1$ to 6) & Iter. on $\Omega_7$ &Iter. on $\Omega$ & Iter. & Order\\
\hline
\multicolumn{7}{|c|}{Case 1: A unit spherical solute region $D_{p}$ containing one central charge only}\\
\hline
0.25	  & $5.77\times 10^{-2}$ & 8.54 $\approx 9$  &7.89 $\approx 8$  & 11.8 $\approx 12$ & 5 & -\\
\hline
$0.25/2$	  & $1.56\times 10^{-2}$ & 8.13 $\approx 8$   & 13.74 $\approx 14$ & 11.4 $\approx 11$ &10 &  1.89\\
\hline
$0.25/4$  & $3.62\times 10^{-3}$ & 8.28 $\approx 8$   &22.98 $\approx 23$  & 11.0 $\approx 11$ & 11 &  2.11\\
\hline
\multicolumn{7}{|c|}{Case 2: The region $D_{p}$ containing 488 point charges from a protein (2LZX).}\\
\hline
0.25	  & $9.10\times 10^{-2}$ & 8.94 $\approx 9$  &7.97 $\approx 8$  & 11.7 $\approx 12$ & 5  & -\\
\hline
$0.25/2$	  & $2.34\times 10^{-2}$ & 8.95 $\approx 9$   & 13.36 $\approx 13$ & 12.2 $\approx 12$ &5 & 1.96 \\
\hline
$0.25/4$  & $5.31\times 10^{-3}$ & 9.61 $\approx 10$   &25.11 $\approx 25$  & 11.0 $\approx 11$ & 5 &  2.14\\
\hline
\end{tabular}}
\end{table}

From Table~\ref{nonlinearV} it can be seen that the errors were reduced almost by three fourths as the mesh size $h$ was decreased by half, indicating that the convergence order of our SMPBE hybrid solver is around 2, well matching the finite element theory. The number of Newton iterations was only up to 11, indicating that our hybrid modified Newton iterative algorithm for computing $\Phit$  retained a fast rate of convergence of the modified Newton method. The average numbers of PCG-MG iterations were about 9 for these three different mesh sizes, numerically confirming that our PCG-MG scheme has a convergence rate independent of the mesh size $h$.
Furthermore, the average numbers of hybrid box iterations were around 11, showing it also attains the convergence rate independent of $h$. 

In these tests, the numbers of PCG-ILU iterations were small, showing the efficiency of PCG-ILU for solving each finite element linear system on the central box $\Omega_{7}$.  Although it increased with the reduction of the mesh size, PCG-ILU was found to take much less CPU runtime than a PCG using an algebraic multigrid preconditioner, called {\em amg\_hypre}, from \texttt{PETSc}. Our test problem sizes might not be large enough to take the advantage of an algebraic multigrid preconditioner.

\subsection{Ionic concentrations for a dipole test model}
\begin{figure}[h]
\centering
\begin{minipage}[b]{0.4\linewidth}
  \includegraphics[width=\linewidth]{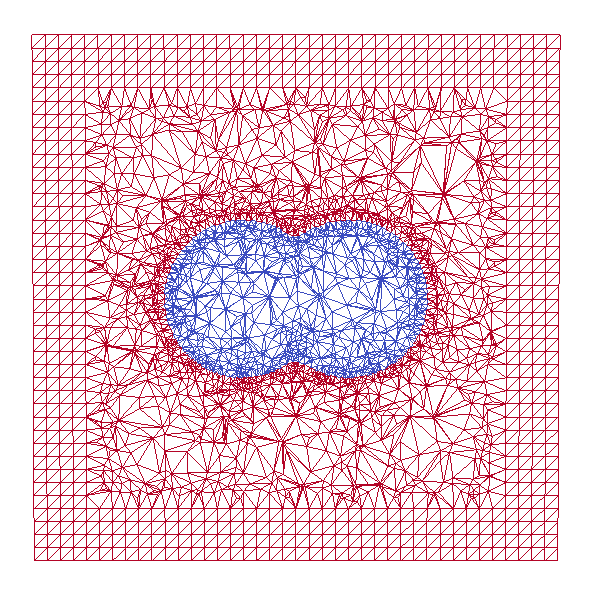}
  \caption{A cross section view on the $xy$ plane of our special mesh of central box $\Omega_{7}$. 
  }
  \label{meshes4d}
\end{minipage}
\quad
\begin{minipage}[b]{0.42\linewidth}
  \includegraphics[width=\linewidth]{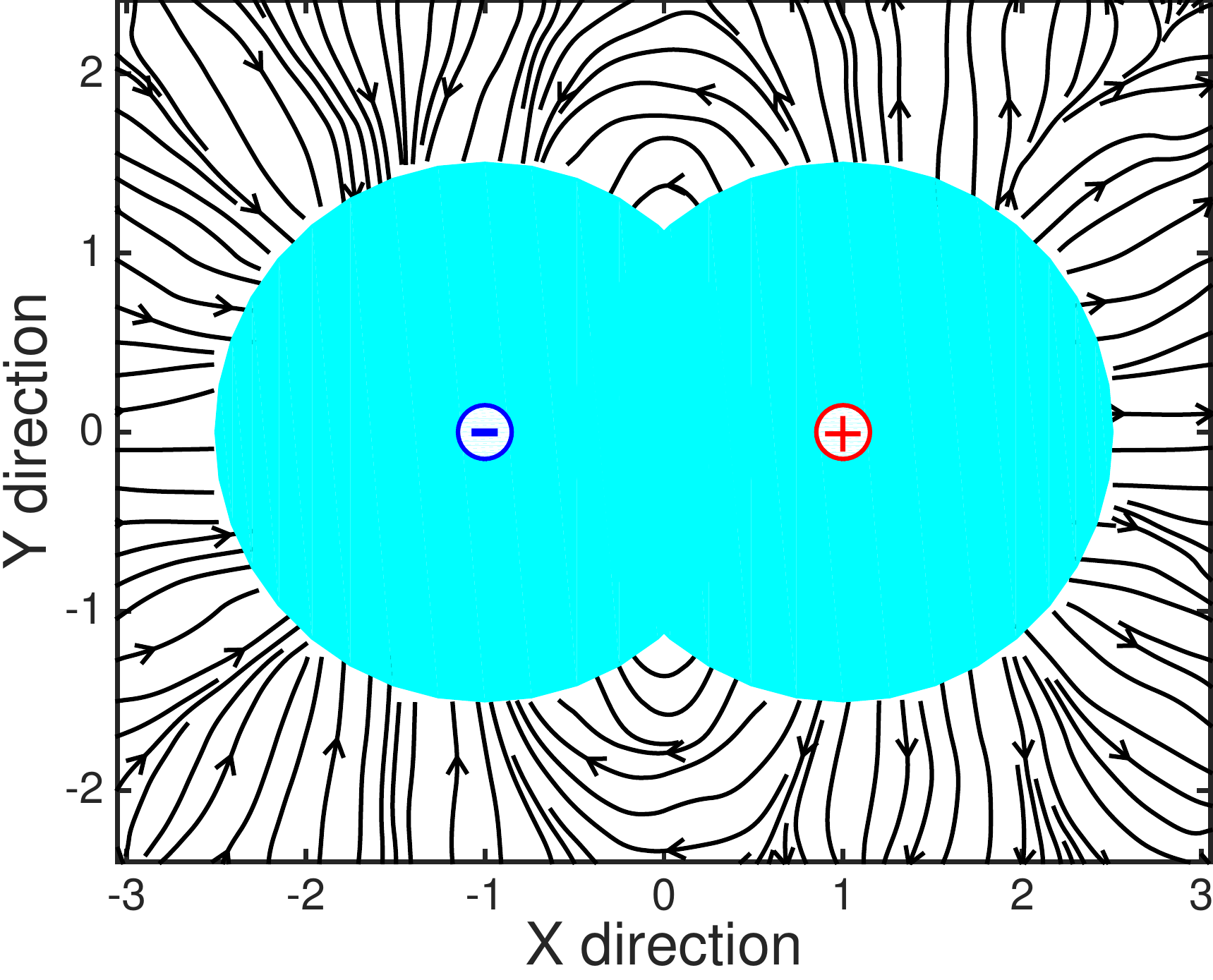}
  \caption{The electrostatic field $E=-\nabla u$ calculated by our SMPBE hybrid solution $u$.   }
  \label{vectorfield}
\end{minipage}
\end{figure}

The SMPBE/PBE test model with $D_{p}$ being a unit ball with a central charge is a common test model to demonstrate that SMPBE can be a better model than PBE in the prediction of ionic concentrations. We did tests on it using the hybrid solver and got the same results as the ones reported in \cite{Li_Xie2014b}. Furthermore, we did tests on a more interesting dipole model, in which $D_p$ consists of two overlapped balls with the same radius $r$ and two opposite central charges. See Figure~\ref{vectorfield} for an illustration.

For a salt solution consisting of sodium (Na$^+$) and chloride (Cl$^-$) ions, the concentrations $C_{Na}$ and $C_{cl}$ of  sodium (Na$^+$) and chloride (Cl$^-$) ions are estimated (in mole per liter) by:
\begin{equation}
\label{c-SMPBE}
    C_{Na} =\left\{\begin{array}{ll} \displaystyle\frac{I_s e^{-u}}{1+2M\Lambda^3\cosh{u}} &  \mbox{ for SMPBE},\\
                   \displaystyle{I_s e^{-u}} &  \mbox{ for PBE},
                   \end{array}\right.   
  C_{cl} = \left\{\begin{array}{ll}  \displaystyle\frac{I_s e^{u}}{1+2M\Lambda^3\cosh{u}} &  \mbox{ for SMPBE},\\
                  \displaystyle{I_s e^{u}}  &  \mbox{ for PBE},
                   \end{array}\right.
\end{equation}
where  $I_s$ is a ionic strength in {\em mole/liter}, and $M$ has been defined in \eqref{alpha-kappa}.

In our numerical tests, we set $I_s=0.1$, $ r=1.5$ \AA, a positive charge of $+3e_c$ at the center $(1,0,0)$, and a negative charge of $-3e_c$ at $(-1,0,0)$. We then constructed an overlapped box partition and meshes with $D=(-4,4)^3$, $\mu=2$, $n=5$ and $m=2$ according to the formulas given in Subsection 4.1 to get $\Omega=(-12,12)^3$ and $\Omega_7=(-5,5)^3$. The meshes of $\Omega$ and $\Omega_{7}$ had  900740 and 56988 mesh points, respectively. A cross-section view of the mesh of $\Omega_7$ on the $xy$-plane is given in Figure \ref{meshes4d}. The boundary value function $g$ was set as zero, which was found numerically to be good enough. 

Figure \ref{vectorfield} displays the electrostatic field $E$ on the $xy$ coordinate plane, which we calculated by the formula $E=-\nabla u$ using a numerical solution $u$ produced by our SMPBE hybrid solver.  From Figure \ref{vectorfield} it can be 
seen that the  electric field lines emanated from the positive charged sphere and extended radially toward the  negative charged sphere. Since the two balls have  an identical quantity of charge, their abilities to alter the space surrounding them are the same. Hence, the number of electric lines are expected to be the same and  the electric lines should occur in a symmetric pattern. As shown in the figure, these basic features of electric field lines were well captured by the numerical solution.  

\begin{figure}[h]
\centering
\begin{subfigure}[b]{0.45\textwidth}
                \centering
                 \includegraphics[width=\textwidth]{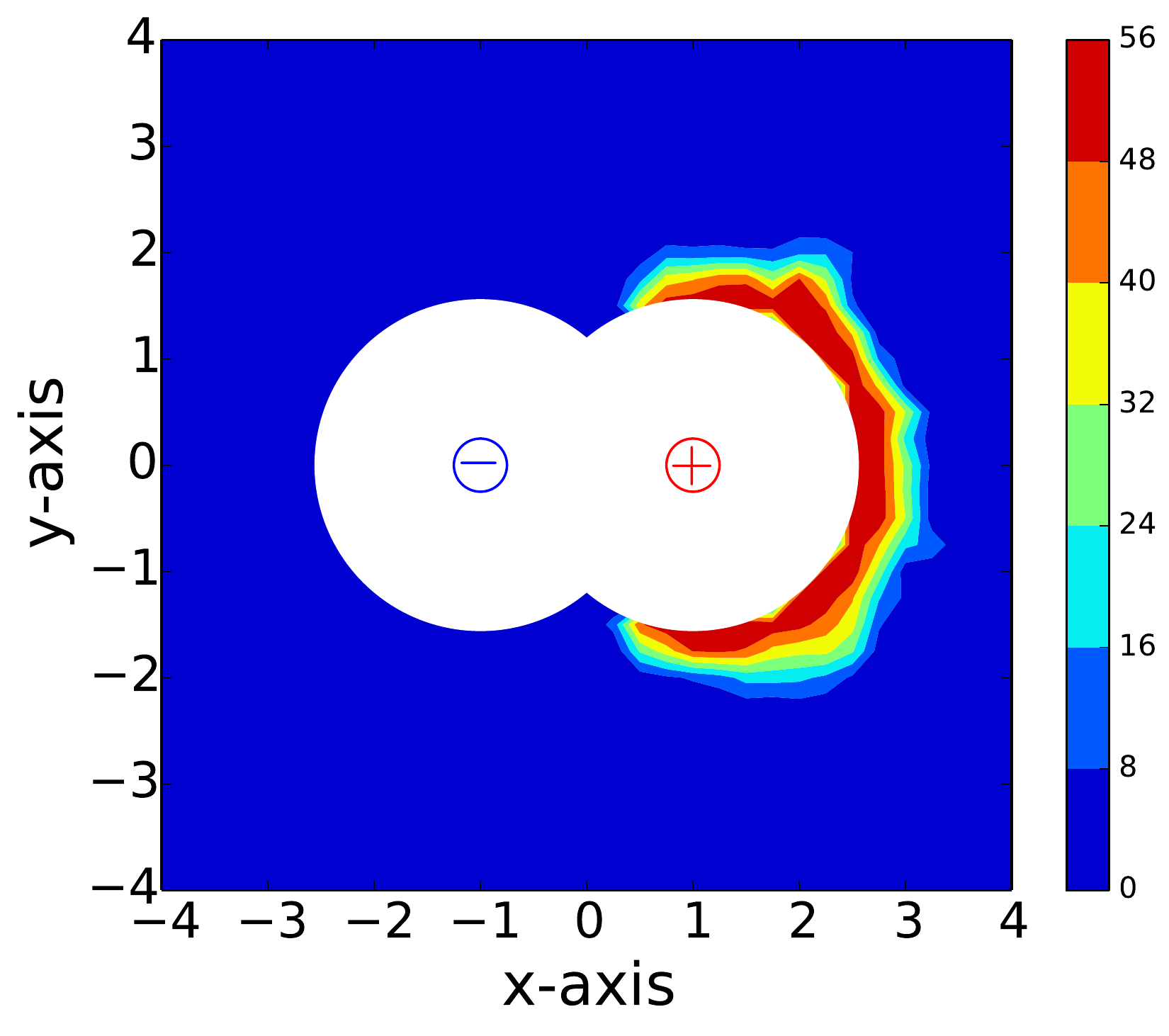}
                \caption{Concentration $C_{cl}$ of anions $Cl^-$}
        \end{subfigure}
        \hspace{2mm}
        \begin{subfigure}[b]{0.45\textwidth}
                \centering
                 \includegraphics[width=\textwidth]{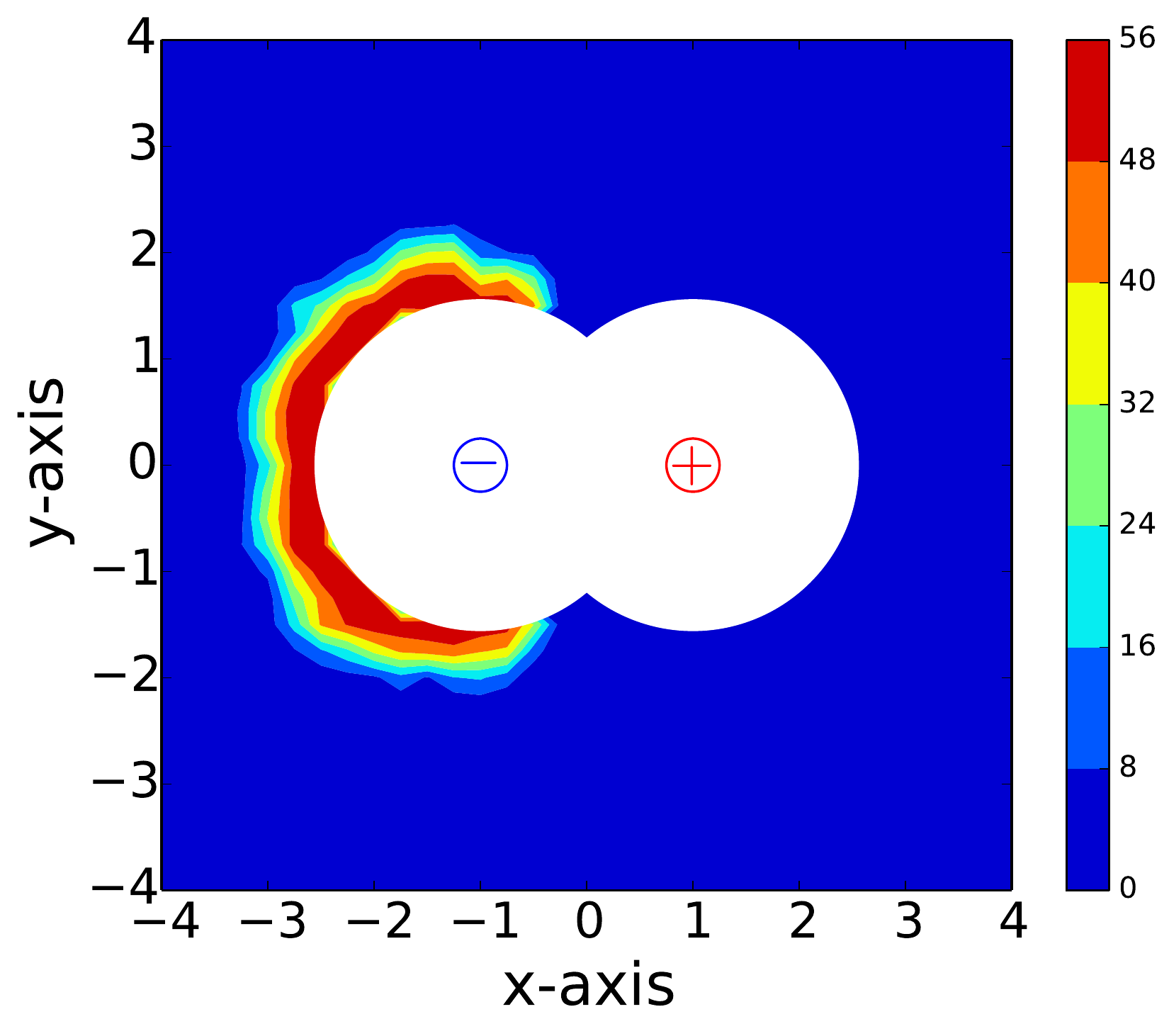}
                \caption{Concentration $C_{Na}$ of cations $Na^+$}
        \end{subfigure}  
\vspace{-0.5em}
\caption{Concentrations of anions ($Cl^-$) and cations ($Na^+$) predicted by our SMPBE hybrid solver for the dipole model on the $xy$-plane.}
\label{concentration_dipole}
\end{figure}

Figure \ref{concentration_dipole} displays the two concentrations $C_{Na}$ and $C_{cl}$ predicted by our SMPBE hybrid solver on the $xy$ coordinate plane. They have reasonably reached a saturation value, 55.2,  claimed in Physics (i.e., $10^{27}/(N_A\Lambda^3)\approx 55.2$).  Interestingly, the predicted values of $C_{Na}$ and 
$C_{cl}$ were distributed symmetrically around the surface of the two balls, well matching the law of electrostatic attraction.

We also repeated the above tests using $\Lambda=0$ to get the results of PBE. The concentrations $C_{Na}$ and $C_{cl}$ were found to be unreasonably large around the spherical surfaces.  

\subsection{Performance tests for proteins}
We made numerical experiments on the  six proteins tested in  \cite{Li_Xie2014b} to show that our SMPBE hybrid program package can improve the performance of the SMPBE finite element program package significantly. These  six proteins have the PDB ID 1CBN, 1SVR, 4PTI, 1AZQ, 1D3X and 1TC3, respectively, which can be downloaded from the Protein Data Bank. Their PDB files were then converted to the PQR files by \texttt{PDB2PQR} \cite{dolinsky2004pdb2pqr}. These six protein molecules have
 642, 1433, 892, 1603, 756, and 2124 atoms,  and 0$e_c$, $-2e_c$, $+6e_c$, $-8e_c$, $-21e_c$, $-35e_c$ net charges, respectively. 
 
In the numerical tests, we used $\mu=4$, $m=2$, and $n=3$ to construct a domain of $\Omega$ and an overlapped box partition of $\Omega$ according to \eqref{domain_construct} for each protein. The over-relaxation parameter $\omega$ was set as 1.225 and 1.015 in solving \eqref{Psi} for $\Psi$ and \eqref{SMPBE_search} for $p_{k}$, respectively. The boundary value function $g$ was set as zero, and the initial iterate $\Phit^{(0)}$ as a numerical solution of the linearized SMPBE reported in \cite{Li_Xie2014b}. 

Because the convergence rate of the solver PCG-MG is independent of the mesh size $h$, our hybrid solver particularly works well on a finer mesh. To demonstrate this feature, we used $m=3$ and $n=4$ to refine the meshes, and then repeated the tests.

We also repeated all the tests by the SMPBE finite element program package. 
Here each tetrahedral mesh of $\Omega$ was produced from the corresponding mesh of each protein for our hybrid solver. That is, each cubic cell of a finite difference mesh was divided into six tetrahedra, making the mesh of the finite element solver to have the same number of mesh points as the mesh used in the hybrid solver. 

The boxes of $D$ produced from our hybrid program for these tests are listed as follows:
\begin{align*} 
D&=(-9.5,28.1)\times(-9.1,28.5)\times(-11.7,25.9) \mbox{ for 1CBN}, \\
D&=(-24.8,27.3)\times(-30.4,21.7)\times(-24.5,27.6)  \mbox{ for 1SVR},\\
D&=(-7.8,38.4)\times(-2.3,43.9)\times(-18.6,27.6)  \mbox{ for 4PTI},\\
D&=(-20.9,40.3)\times(-19.6,41.6)\times(-19.2,42.0)  \mbox{ for 1AZQ},\\
D&=(-21.1,20.7)\times(-21.6,20.2)\times(-20.4,21.4) \mbox{ for 1D3X}, \\
D&=(-30.3,51.0)\times(91.4,172.7)\times(-8.1,73.2)  \mbox{ for 1TC3}.
\end{align*}

\begin{table}[t]
\caption{A comparison of the performance in CPU time (in seconds) of our hybrid solver (Hybrid) with that of the  finite element solver (FE) reported  in \cite{Li_Xie2014b} in the calculation of component functions $G$, $\Psi$ and $\Phit$ of SMPBE solution $u$, including the total CPU time (exclusion of the time for finite element mesh generation). }
\label{protein-case}
\begin{tabular}{|c|c|c|c|c|c|c|c|c|}
 \hline
Number of   & \multicolumn{2}{|c|}{Find $G$ \& $\nabla G$ } & \multicolumn{2}{|c|}{Find $\Psi$} & \multicolumn{2}{|c|}{Find $\Phit$}   & \multicolumn{2}{|c|}{Total Time}    \\ \cline{2-9}
mesh points  &   Hybrid & FE &  Hybrid & FE &  Hybrid  & FE &  Hybrid & FE \\ 
 \hline
  \multicolumn{9}{|c|}{ Protein with PDB ID 1CBN (642 atoms and 0$e_{c}$ net charge)} \\ \hline
  77,969   & 0.42 & 0.76 & 0.91 & 1.99 & 2.32 & 8.02  & 4.17 & 11.33 \\
\hline
 537,953   & 2.63 & 5.23 & 4.71 & 14.97 & 8.31 & 50.21  & 17.81 & 74.09 \\
\hline
 \multicolumn{9}{|c|}{ Protein with PDB ID 1SVR (1433 atoms and $-2e_{c}$ net charge)} \\ \hline
 89,850  &1.40   &2.37   &1.77 &2.53 &5.05 &10.71   &9.14 & 16.25 \\
\hline
550,170  &6.96   &14.49   &5.96 &15.46 &11.40 &60.10   &27.02 & 93.84 \\
\hline
 \multicolumn{9}{|c|}{ Protein with PDB ID 4PTI (892 atoms and $+6e_{c}$ net charge)} \\ \hline
  81,356  &  0.74  & 1.33    & 1.20  &2.12 & 3.35  & 9.47  & 5.92 & 13.54 \\
\hline
541,329  &4.19    & 8.84    & 5.47  &14.96& 10.21  & 50.45  & 22.30 & 77.97 \\
\hline
 \multicolumn{9}{|c|}{ Protein with PDB ID 1AZQ (1603 atoms and $-8e_{c}$ net charge)} \\ \hline
 89,089 &1.54   &2.62  &1.69  &2.40 &4.73  &10.54 &8.87  &16.20 \\
\hline
549,279 &7.75   &16.14  &5.94  &15.27 & 13.21 &57.29 &29.76  &92.49 \\
\hline
 \multicolumn{9}{|c|}{ Protein with PDB ID 1D3X (756 atoms and $-21e_{c}$ net charge)} \\ \hline
 82,897 &  0.64 &1.15  &1.32  &2.15 &4.14  &10.65 &6.77  &14.54 \\
\hline
542,878 &3.58   & 7.54  &5.86  &14.76 &12.78  &57.48 &24.59  &83.53 \\
\hline
 \multicolumn{9}{|c|}{ Protein with PDB ID 1TC3 (2124 atoms and $-35e_{c}$ net charge)} \\ \hline
 101,944 &2.54   & 3.99  &2.90  &2.91 &10.05  &16.81  &16.90  &24.45  \\
\hline
564,871 &10.88   &21.98  &7.58  &15.65 &23.50  &85.29  &45.45  &126.82  \\
\hline
 \end{tabular}
\end{table}
 
Table \ref{protein-case} reports these numerical test results. From it we can see that our hybrid solver reduced the total CPU runtime of the SMPBE finite element solver  from $31\%$ to $63\%$ on the coarse meshes and $64\%$ to $76\%$ on the fine meshes, showing that our hybrid solver can significantly improve the performance of the original finite element solver, especially in the case  of fine meshes. Besides, the memory usage of the finite element solver can be reduced sharply by our hybrid solver, since our PCG-MG solver does not require any array to store the mesh or coefficient data of a finite difference linear system.

\subsection{Electrostatic solvation free energy calculations}
One important application of SMPBE is to predict the electrostatic solvation free energy $\triangle E$ for a protein in an ionic solvent, which describes the energy change for the protein to move from the reference state to  the solvated state. By the solution decomposition \eqref{solutionSplit},  $\triangle E$ can be estimated in kilocalorie per mole  (kcal/mol) by the following formula
\begin{equation}
\label{Solv}
\triangle E = \frac{N_A}{4184}\cdot\frac{k_B T}{2}\sum_{j=1}^{n_p} z_j\left(\Psi(\rr_j)+\Phit(\rr_j)\right).
\end{equation}

Clearly,  numerical solutions $\Psi_{h}$ and $\Phit_{h}$ give a numerical value $\triangle E_{h}$, and $\triangle E_{h}\to \triangle E$ if $\Psi_{h} \to \Psi$ and
$\Phit_{h}\to \Phit$ as $h\to 0$. Hence, a numerical calculation of $\triangle E$ provides us with a way to check the convergence behavior of our hybrid solver.

\begin{figure}[h]
\centering
\begin{subfigure}[b]{0.448\textwidth}
                \centering
                 \includegraphics[width=\textwidth]{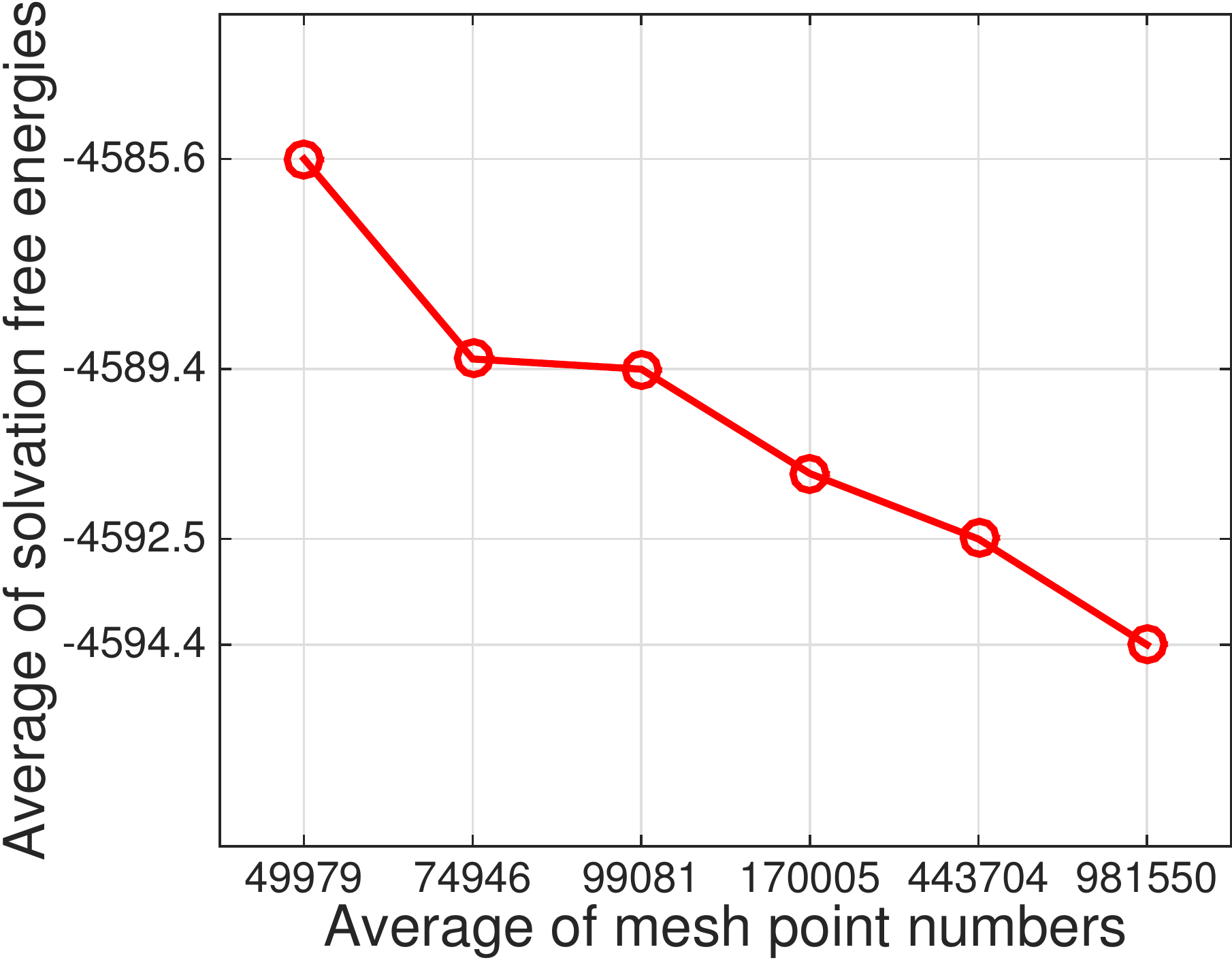}
        \end{subfigure}
        \hspace{2mm}
        \begin{subfigure}[b]{0.45\textwidth}
                \centering
                 \includegraphics[width=\textwidth]{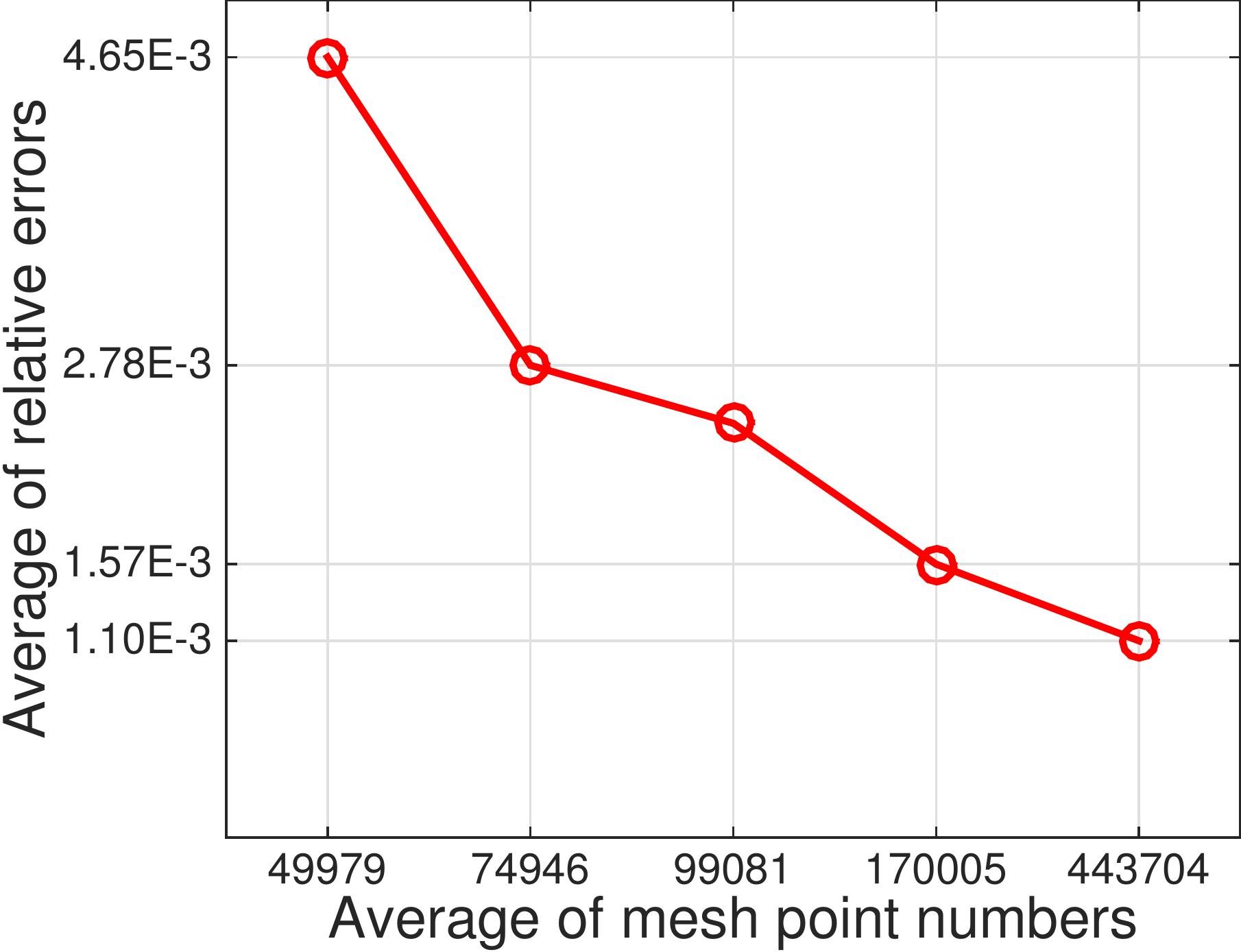}
        \end{subfigure}  
\vspace{-0.5em}
\caption{Numerical behavior of the new SMPBE hybrid solver on different meshes in the calculation of electrostatic solvation free energy for the 216 biomolecules.}
\label{SolVConvergence}
\end{figure}

To do so, we calculated $\triangle E$ by the formula  \eqref{Solv} for a set of 216 biomolecules with atom numbers from 506 to 69711 (including proteins, protein-protein complex, and nucleic acid) downloaded from  Prof. Ray Luo's website {\em http://rayl0.bio.uci.edu/rayl/}. Here, each domain $\Omega$ was selected to be 2 times larger than box D (i.e., $\mu=2$), and the zero boundary condition was used.
For each biomolecule, we calculated $\triangle E$ using six successively refined meshes. The averages of the mesh point numbers over the 216 meshes of whole domain $\Omega$ were found to be 49979, 74946, 99081, 170005, 443704, and 981550 for the six sets of meshes, respectively. Since the analytical value of $\triangle E$ is unknown,  we took the numerical value of $\triangle E$ calculated from the finest mesh as the reference to calculate an approximate relative error of $\triangle E_{h}$. To simplify the display of numerical results, we calculated the average of 216 relative errors for each set of meshes as reported in Figure \ref{SolVConvergence}. From Figure \ref{SolVConvergence} we can see that the relative errors were changed less than 0.004 only, implying that our SMPBE hybrid solver has good properties of the numerical stability and convergence
in the calculation of electrostatic solvation free energy.

\subsection{Binding free energy calculations}
Similar to the classic PBE \cite{Honig2007PBE,sitkoff1994accurate}, SMPBE can also be used to study the salt dependence of the binding free energy for a complex molecule.   
For a complex $C$ consisting of molecules $A$ and $B$, the binding free energy $E_b(I_s)$ is defined by 
\[E_b(I_s) = E(C,I_s)-E(A,I_s)-E(B,I_s),\]
where  $E(X,I_s)$ denotes an electrostatic free energy of molecule $X$ in a solvent with the ionic strength $I_s$. 
From the counterion condensation theory \cite{fenley2010revisiting,manning1978molecular}, it is known that $E_b$ can be transformed by the  variable change, $\xi=\ln I_s$, to  a linear function of $\xi$ as follows
\begin{equation}
\label{bindingslope}
E_b = m \xi + b,
\end{equation}
where $m$ and $b$ are constants to be determined. In \cite{breslauer1987enthalpy},  a scaled slope $m_s$ is defined by
\[  m_s = - m /(N_Ak_BT),\]

We made tests on  a DNA-drug complex represented in PDB ID 1D86 using the PQR files from \cite{fenley2010revisiting}. A chemical experimental value, $-1.51$,  of $m_s$ was given in  \cite[Table 3]{breslauer1987enthalpy}.
In the numerical tests, the electrostatic free energy $E(X,I_{s})$ was calculated by our hybrid solver according to the formula of \eqref{Solv}. 

To produce a numerical prediction to the scaled slope $m_{s}$, we calculated the binding free energy $E_b$ using  the following 11 different values of $I_{s}$:
\[  I_{s,j}=e^{\xi_{j}} \quad  \mbox{ with }  \quad \xi_{j}=-3+0.2 j \quad \mbox{ for }   j=0,1,2,\ldots,9,10,\]
and then generate a best-fit line by a linear regression program (downloaded from the APBS website) to yield a predicted value of $m_{s}$.

\begin{figure}[t]
\centering
\begin{minipage}[b]{0.45\linewidth}
  \includegraphics[width=1.0\linewidth]{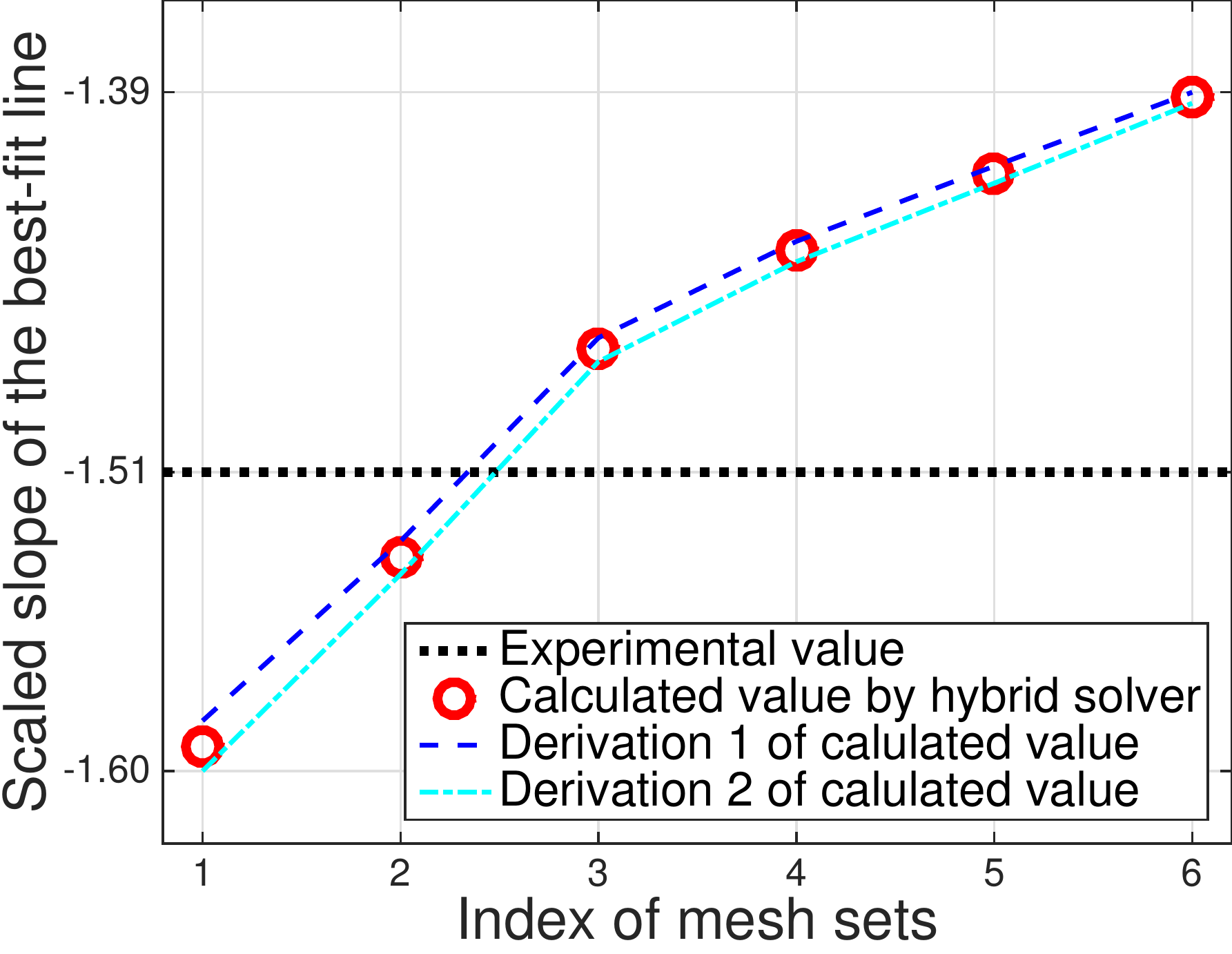}
  \caption{The scaled slops $m_{s}$ of the six best-fit lines calculated by our SMPBE hybrid solver  for a DNA-drug complex represented in PDB ID 1D86 based on the six sets of meshes listed in Table~\ref{table-meshes}.  }
  \label{slop1}
\end{minipage}%
\qquad
\begin{minipage}[b]{0.45\linewidth}
  \includegraphics[width=0.934\linewidth]{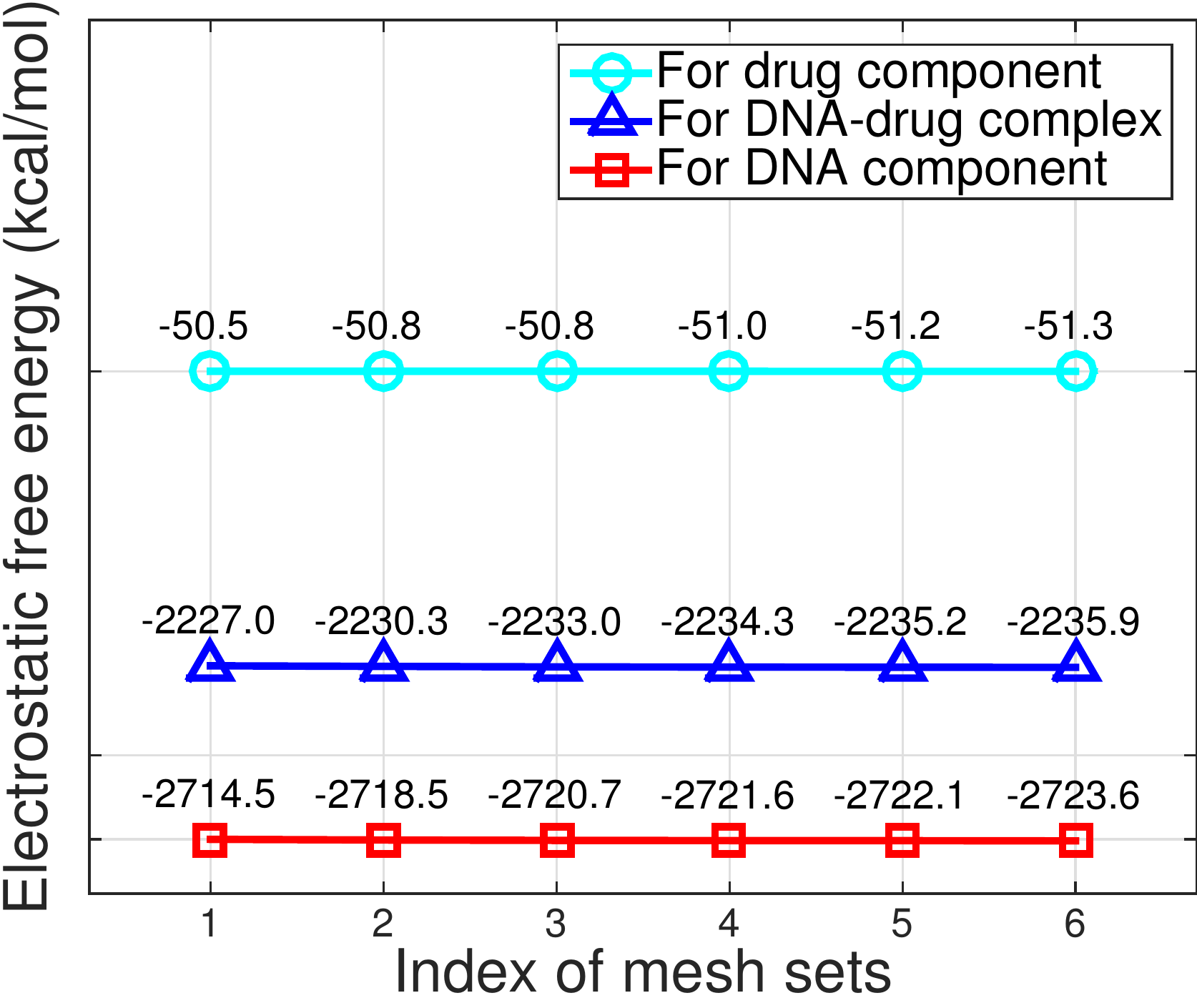}
  \caption{Electrostatic free energies  (numbers on the top of each curve)  of the complex, DNA, and drug calculated by the hybrid solver on the same six mesh sets of Fig.~\ref{slop1} using $I_s \approx 0.11$. }
  \label{slop2}
\end{minipage}
\end{figure}

Moreover, we constructed six different sets of meshes and repeated the above calculation in order to study the sensitivity of our hybrid solver to discretization error. Here each set contains three meshes:  the first one for the complex 1D86, the second one for its DNA component, and the third one for its drug component. In each mesh, a domain of $\Omega$ was set to be four times larger than  a domain of $D$, i.e., $\mu=4$. The interface  $\Gamma$ produced from the first mesh were shared by the other five meshes to ensure that the same interface problem was solved on each mesh. The numbers of mesh points of these meshes are listed in Table \ref{table-meshes}. Numerical results are displayed in Figures~\ref{slop1} and \ref{slop2}.

\begin{table}[ht]
\centering
\begin{tabular}{|c|c|c|c|}
 \hline
\text{Mesh set}& \multicolumn{3}{|c|}{  Number of mesh points}  \\ 
 \cline{2-4}
 index &  Mesh for complex  & Mesh for DNA   & Mesh for drug  \\
\hline
1&  85195  & 86435 & 73051   \\
 \hline
2& 154705 & 155903 & 137797\\
 \hline
3& 259903 & 261512 & 233983   \\
 \hline
4& 411583 & 413961 & 367030 \\
 \hline
5& 613983 & 616328 & 542747 \\
 \hline
 6& 4252682 & 4255163 & 4163790\\
 \hline
\end{tabular}
\caption{Mesh data for the six sets of meshes used in the calculation of binding free energy for a DNA-drug complex represented in PDB ID 1D86. }
\label{table-meshes}
\end{table}

From Figure \ref{slop1} it can be seen that the numerical values of the scaled slop $m_{s}$ increased from $-1.6$ to $-1.39$ across the experimental value of $-1.51$, as the mesh set index was changed from 1 to 6 (or a mesh became finer and finer). The derivations from the best-fit line calculation were all very small. These tests showed that our hybrid solver behaved stably in the numerical calculation of binding free energy. In Figure \ref{slop2}, we displayed the electrostatic free energy values of $E(X,I_{s})$ calculated by our hybrid solver using $I_{s} = e^{-2.2} \approx 0.11$, which showed a tendency of convergence for $X$ being the DNA-drug complex, DNA component, or drug component.
This important property makes it easy for us to produce a good predicted value of $m_{s}$ on a properly constructed mesh set. 
\section*{Acknowledgements}
This work was partially supported by the National Science Foundation, USA, through grant DMS-1226259. 
The authors would like to thank Dr. Marcia O. Fenley for the PQR files of the DNA-drug complex represented in PDB ID 1D86. 


\end{document}